\documentclass[12pt,twoside]{article}
\usepackage{amsmath,amssymb,amsfonts,amstext,amsthm,amsxtra}
\usepackage{xypic}
\usepackage{color}
\usepackage[dvips]{epsfig,graphics}
\usepackage[all,knot,arc,poly]{xy}
\usepackage[english]{babel}
\usepackage{graphics}

\date{}

\begin{document}

\centerline{}

\centerline{}

\centerline {\Large{\bf First Non-Abelian Cohomology of Topological Groups}}

\centerline{}
\newcommand{\mvec}[1]{\mbox{\bfseries\itshape #1}}

\centerline{\bf {H. Sahleh$^1$ and H. E. Koshkoshi$^2$}}

\centerline{}

\centerline{$^1$$^,$$^2$Department of  Mathematics, Faculty of Mathematical Sciences, University of Guilan}
\centerline{P. O. Box 1914, Rasht-Iran}
\centerline{$^1$E-mail: sahleh@guilan.ac.ir}
\centerline{$^2$E-mail: h.e.koshkoshi@guilan.ac.ir}

\centerline{}

\newtheorem{Theorem}{\quad Theorem}[section]

\newtheorem{Definition}[Theorem]{\quad Definition}

\newtheorem{Corollary}[Theorem]{\quad Corollary}
\newtheorem{Proposition}[Theorem]{\quad Proposition}
\newtheorem{Lemma}[Theorem]{\quad Lemma}
\newtheorem{Remark}[Theorem]{\quad Remark}
\newtheorem{Example}[Theorem]{\quad Example}
\newtheorem{Notice}[Theorem]{\quad Notice}
\centerline{}
\centerline{\bf Abstract}
{\emph{Let $G$ be a topological group and $A$ a topological $G$-module (not necessarily abelian). In this paper, we define $H^{0}(G,A)$ and $H^{1}(G,A)$ and will  find a six terms exact cohomology sequence involving $H^{0}$ and $H^{1}$. We will extend it to a seven terms exact sequence of cohomology up to dimension two. We find a criterion  such that vanishing of  $H^{1}(G,A)$ implies the connectivity of $G$. We show that if $H^{1}(G,A)=1$, then all complements  of $A$ in the semidirect product $G\ltimes A$ are conjugate. Also as a result, we prove that if $G$ is a compact Hausdorff group and $A$ is a locally compact almost connected Hausdorff group with the trivial maximal compact subgroup then, $H^{1}(G,A)=1$.}}

\par\textbf{2000 MSC No:} Primary 22A05, 20J06; Secondary 18G50.

\section{Introduction}
Let $G$  and $A$ be topological groups. It is said that $A$ is a topological $G$-module, whenever $G$  continuously acts  on the  left of $A$. For all $g\in G$ and $a\in A$ we denote the action of $g$ on $a$ by $^{g}a$.

\par In section 2, We define $H^{0}(G,A)$ and $H^{1}(G,A)$.
\par In section 3, we define  the covariant functor  $H^{i}(G,-)$ for $i=0,1$ from the category of topological $G$-modules to the category of pointed sets. Also, we define two connecting maps $\delta^{0}$ and $\delta^{1}$.
\par A classical result of Serre [6], asserts that if $G$ is a topological group and $1\rightarrow A\rightarrow B\rightarrow C\rightarrow 1$ a central short exact sequence of discrete $G$-modules then, the sequence
$1\rightarrow H^{0}(G,A)\rightarrow H^{0}(G,B)\rightarrow H^{0}(G,C)\rightarrow H^{1}(G,A)\rightarrow H^{1}(G,B)\rightarrow H^{1}(G,C)\rightarrow  H^{2}(G,A)$ is exact.
\par In section  4, we generalize the above result to
the case of  arbitrary topological $G$-modules (not necessarily discrete).
\par We show that if $G$ is a connected group and $A$  a totally disconnected group then, $H^{1}(G,A)=1$.
\par In section 5,  we show that if $G$ has an open component (for example $G$ with the finite number of components) and for every discrete (abelian) $G$-module $A$ $H^{1}(G,A)=1$ then, $G$ is a connected group.
\par In section 6, we show that vanishing of $H^{1}(G,A)$ implies that the complements of $A$ in the (topological) semidirect product
 $G\ltimes A$, are conjugate.
\par In section 7, we prove that, if $G$ is a compact Hausdorff group and  $A$ a locally compact almost connected Hausdorff group
then there exists a  $G$-invariant maximal compact subgroup $K$ of $A$  such that  the natural map
 $\iota_{1}^{*}:H^{1}(G,K)\rightarrow H^{1}(G,A)$ is onto. As a result, if $G$ is compact Hausdorff and $A$ is a locally compact almost connected Hausdorff group with  trivial maximal compact subgroup then, $H^{1}(G,A)=1$.

\par All topological groups are arbitrary (not necessarily abelian).
We assume that $G$ acts on itself by conjugation.
The center of a group $G$  and the set of all continuous homomorphisms of $G$ into $A$ are denoted by $Z(G)$ and $Hom_{c}(G,A)$,
respectively. The topological isomorphism is denoted by $"\simeq" $.
\par Suppose that $A$ is an abelian topological $G$-module.\\
Take $\tilde{C}^{0}(G,A)=A$ and  for every positive integer $n$, let $\tilde{C}^{n}(G,A)$ be  the set of continuous maps $f:G^{n}\rightarrow A$ with the coboundary map $\tilde{\delta}^{n}:\tilde{C}^{n}(G,A)\rightarrow \tilde{C}^{n+1}(G,A)$  given by
\begin{center}
$\tilde{\delta}^{n}f(g_{1},...,g_{n+1})=\ ^{g_{1}}f(g_{2},...,g_{n+1})+\sum (-1)^{i}f(g_{1},...,\underbrace{g_{i}g_{i+1}}_{ith} ,...,g_{n+1})+(-1)^{n+1}f(g_{1},...,g_{n}).$

\end{center}
The $n$th cohomology of $G$ with coefficients in $A$ in the sense of Hu [5], is the abelian group
$$H^{n}(G,A)=Ker\tilde{\delta}^{n}/Im\tilde{\delta}^{n-1}.$$

\section{ $H^{0}(G,A)$ and $H^{1}(G,A)$}
 Let $G$ be a topological group and $A$ a topological $G$-module.
\begin{Definition}We define $H^{0}(G,A)=\{a|a\in A, ^{g}a=a, \forall g\in G\}$, i.e., $H^{0}(G,A)=A^{G}$, the set of $G$-fixed elements of $A$.
\end{Definition}
\begin{Definition}  A  map $\alpha:G\rightarrow A$ is called a continuous derivation  if $\alpha$ is continuous and $$\alpha(gh)=\alpha(g)^{g}\alpha(h),  \forall g,h\in G.$$
\end{Definition}
The set of all continuous derivations from $G$ into $A$ is denoted by $Der_{c}(G,A)$. Two continuous derivations  $\alpha, \beta$ are cohomologous, denoted by $\alpha \sim \beta$, if there is  $a\in A$ such that
\begin{center}
$\beta(g)=a^{-1}\alpha(g)^{g}a$,  for all $g\in G$.
\end{center}
\par It is easy to show that $\sim$ is an equivalence relation.
Now we define $$H^{1}(G,A)=Der_{c}(G,A)/\sim.$$
\begin{Notice}There exists the trivial continuous derivation  $\alpha_{0}:G\rightarrow A$ where $\alpha_{0}(g)=1$; Hence, $H^{1}(G,A)$ is nonempty.
In  general, $H^{1}(G,A)$ is not a group. Thus, we will view $H^{1}(G,A)$ as a pointed set with the basepoint $\alpha_{0}$.
\par Note that $H^{0}(G,A)$ is a subgroup of $A$, so it is a pointed set with the basepoint 1.  Also if $A$ is a Hausdorff group,
 then, $H^{0}(G,A)$ is a closed subgroup of $A$.
\end{Notice}
%
%
%
%
%
%
\begin{Remark}
(i) If $A$ is an abelian group then, $H^{1}(G,A)$ is the first (abelian) group cohomology in the sense of Hu, i.e., it is the group of all continuous derivations of $G$ into $A$ reduced modulo the inner derivations. [5]
\\ \\(ii) If $A$ is a trivial topological $G$-module then, $H^{1}(G,A)=Hom_{c}(G,A)/\sim$. Here $\alpha \sim \beta$ if $\exists \ a\in A$ such that $\beta(g)=a^{-1}\alpha(g)a, \forall g\in G$.
\\\\(iii) Let $G$ be a connected group and $A$ a totally disconnected group then, $H^{1}(G,A)=1$.
\end{Remark}
{\bf Proof.}  (i) and (ii) are obtained  from the definition of $H^{1}(G,A)$.
 \\(iii): If $\alpha \in Der_{c}(G,A)$ then, $\alpha(1)=1$. On the other hand $G$ is a connected group and $A$ is totally disconnected.  So, $\alpha=\alpha_{0}$. Thus,  $H^{1}(G,A)=1$.

\section{ $H^{i}(G,-)$ as a Functor and the Connecting Map $\delta^{i}$ for $i=0,1$}
In this section we define two covariant functors $H^{0}(G,-)$ and $H^{1}(G,-)$ from   the category of topological $G$-modules $_{G}\mathcal{M}$ to    the category of pointed sets $\mathcal{PS}$. Furthermore, We will define the connecting maps $\delta^{0}$ and $\delta^{1}$.
\\
\par Let $A, B$ be topological $G$-modules and $f:A\rightarrow B$ a continuous $G$-homomorphism. We define  $H^{i}(G,f)=f_{i}^{*}:H^{i}(G,A)\rightarrow H^{i}(G,B)$, $i=0,1$, as follows:
\par For $i=0$,  take $f_{0}^{*}=f|_{A^{G}}$. This gives a homomorphism from $H^{0}(G,A)$ to $H^{0}(G,B)$, since $f$ is a homomorphism of $G$-modules .
 So if $a \in A^{G}$, then, $^{g}f(a)=f(^{g}a)=f(a)$, for each $g\in G$. Hence,
 $f(a)\in B^{G}$, i.e., $f_{0}^{*}$ is  well-defined.
\par For $i=1$, we define $f_{1}^{*}$ as follows:
\\
For simplicity, we write $\alpha$ instead of $[\alpha] \in H^{1}(G,A)$.
 \\If $\alpha \in H^{1}(G,A)$, then, take $f_{1}^{*}(\alpha)=f\circ\alpha$. Now if $g, h\in G$, then,
 $$f_{1}^{*}(\alpha)(gh)=f(\alpha(gh))=f(\alpha(g)^{g}\alpha(h))=f(\alpha(g))f(^{g}\alpha(h))=f_{1}^{*}(\alpha)(g)^{g}f_{1}^{*}(\alpha)(h).$$ So, $f_{1}^{*}(\alpha)$ is a continuous derivation.
\par Moreover, if $\alpha, \beta\in H^{1}(G,A)$ are cohomologous then, there is $a\in A$ such that
$\beta(g)=a^{-1}\alpha(g) ^{g}a$. Hence, $f(\beta(g))=f(a)^{-1}f(\alpha(g))^{g}f(a)$. So, $f_{1}^{*}(\alpha) \sim f_{1}^{*}(\beta)$.
\\The fact that  $H^{i}(G,-)$ is a functor  follows from the definition of $f_{i}^{*}$, ($i=0, 1$).  Also  $H^{0}(G,-)$ is a covariant functor from $_{G}\mathcal{M}$ to the category of topological groups $\mathcal{TG}$.
\\
\par Suppose that $ \xymatrix{ 1 \ar @{->}[r] &A \ar @{->} [r]^{\iota}
 &B\ar @{->} [r]^{\pi}& C \ar @{->}[r]& 1}$ is an exact sequence of topological $G$-modules and continuous $G$-homomorphisms such that $\iota$ is an embedding. Thus, we can identify $A$ with $\iota(A)$.
  \par
 Now we define a coboundary map $\delta^{0}:H^{0}(G,C)\rightarrow H^{1}(G,A)$.
 \\Let $c\in H^{0}(G,C)$,  $b\in B$ with $\pi(b)=c$. Then, we define $\delta^{0}(c)$ by $\delta^{0}(c)(g)=b^{-1}$$^{g}b, \forall g\in G$.
 It is obvious that $\delta^{0}(c)$ is a continuous derivation. Let $b^{'} \in B$, $\pi(b^{'})=c$. Then, $b^{'}=ba$ for some $a\in A$. So,

\begin{center}
$(b^{'})^{-1}$$^{g}b^{'}=a^{-1}b^{-1}$$^{g}b^{g}a=a^{-1}\delta^{0}(c)(g)^{g}a$.
\end{center}
Thus, the derivation obtained from $b{'}$ is cohomologous in $A$ to the one obtained from $b$,
i.e., $\delta^{0}$ is well-defined.
\par Now,  suppose that $ \xymatrix{ 0 \ar @{->}[r] &A \ar @{->} [r]^{\iota}
 &B\ar @{->} [r]^{\pi}& C \ar @{->}[r]& 1}$  is a central exact sequence of $G$-modules and continuous $G$-homomorphisms such that $\iota$ is a homeomorphic embedding and in addition $\pi$ has a continuous section $s:C\rightarrow B$, i.e., $\pi s=Id_{C}$.
 \par We construct a coboundary map $\xymatrix{ H^{1}(G,C)\ar @{->}[r]^{\delta^{1}}& H^{2}(G,A)}$. Here $H^{2}(G,A)$ is defined in the sense of Hu [5]. By assumption $\iota(A) \subset Z(B)$, so, $A$ is an abelian topological $G$-module.
 \\ Let $\alpha \in H^{1}(G,C)$ and $s:C\rightarrow B$ be a continuous section for $\pi$. Define $ \delta^{1}(\alpha)$ via $\delta^{1}(\alpha)(g,h)=s\alpha(g)$ $^{g}(s\alpha(h))(s\alpha(gh))^{-1}$. It is clear that $ \delta^{1}(\alpha)$ is a continuous map.
 \par We show that $\delta^{1}(\alpha)$ is a factor set with values in $A$,
   and independent of the choice of the continuous section $s$. Also $\delta^{1}$ is well-defined.
 \\Since $\alpha$ is a derivation, we have:
 \begin{center}
 $\pi(\delta^{1}(\alpha)(g,h))=\pi($$s\alpha(g)$ $^{g}(s\alpha(h))(s\alpha(gh))^{-1})=\alpha(g)^{g}\alpha(h)(\alpha(gh))^{-1}=1$.
 \end{center}
 Thus, $\delta^{1}(\alpha)$ has values in $A$.
 \par Next, we show that $\delta^{1}(\alpha)$ is a factor set, i.e.,
$$^{g}\delta^{1}(\alpha)(h,k)\delta^{1}(\alpha)(g,hk)=\delta^{1}(\alpha)(gh,k)\delta^{1}(\alpha)(g,h), \ \forall g, h, k \in G.\eqno{(3.1)}$$

 First we calculate the left hand side of  (3.1). For simplicity, take $b_{g}=s\alpha(g)$, $\forall g\in G$. Since $A\subset Z(B)$, thus,

\begin{center}
 $^{g}\delta^{1}(\alpha)(h,k)\delta^{1}(\alpha)(g,hk)=\-^{g}(b_{h}\-^{h}b_{k}b_{hk}^{-1})(b_{g}\-^{g}b_{hk}b_{ghk}^{-1})=
 b_{g}\-^{g}(b_{h}\-^{h}b_{k}b_{hk}^{-1})\-^{g}b_{hk}b_{ghk}^{-1}
$ \end{center}
\begin{center}
$=b_{g}\-^{g}(b_{h}\-^{h}b_{k})^{g}b_{hk}b_{ghk}^{-1}=b_{g}\-^{g}b_{h}\-^{gh}b_{k}b_{ghk}^{-1}$,
\end{center}
On the other hand,
\begin{center}
$\delta^{1}(\alpha)(gh,k)\delta^{1}(\alpha)(g,h)=(b_{gh}\-^{gh}b_{k}b_{ghk}^{-1})(b_{g}\-^{g}b_{h}b_{gh}^{-1})
=b_{g}\-^{g}b_{h}\-^{gh}b_{k}b_{ghk}^{-1}$.
\end{center}
Therefore, $\delta^{1}(\alpha)$ is a factor set.
\par Next, we prove that $\delta^{1}(\alpha)$ is independent of the choice of the continuous section. Suppose that $s$ and $u$ are continuous sections for $\pi$. Take $b_{g}=s\alpha(g)$ and $b^{'}_{g}=u\alpha(g)$, for a fixed $\alpha \in Der_{c}(G,C)$. Since $\pi(b_{g}^{'})=\alpha(g)=\pi(b_{g})$, then, $b_{g}^{'}=b_{g}a_{g}$ for some $a_{g}\in A$. Obviously the function $\kappa:G\rightarrow A$, defined by $\kappa(g)=a_{g}$, is continuous. Thus, \begin{center}
    $(\delta^{1})^{'}(\alpha)(g,h)=b^{'}_{g}\-^{g}b^{'}_{h}b^{'}_{gh}=b_{g}\kappa(g)\-^{g}b_{h}\-^{g}\kappa(h)(\kappa(gh))^{-1}b^{-1}_{gh})$
\end{center}
\begin{center}
    =$(\kappa(g)^{g}\kappa(h)(\kappa(gh))^{-1})(b_{g}\-^{g}b_{h}b_{gh}^{-1})=\tilde{\delta}^{1}(\kappa)(g,h)\delta^{1}(\alpha)(g,h),$
\end{center}
where $\tilde{\delta}^{1}(\kappa)(g,h)=\-^{g}\kappa(h)(\kappa(gh))^{-1}\kappa(g).$ \\The coboundary map $\tilde{\delta}^{1}:\tilde{C}^{1}(G,A)\rightarrow \tilde{C}^{2}(G,A)$ is defined as in [5].
 Consequently, $\delta^{1}(\kappa)$ and $(\delta^{1})^{'}(\kappa)$ are cohomologous.
\\
\par Suppose that $\alpha$ and $\beta$ are cohomlologous in $Der_{c}(G,A)$. Then, there is  $c\in C$ such that $\beta(g)=c^{-1}\alpha(g)^{g}c$, $\forall g\in G$.
\\Let $s:C\rightarrow A$ be a continuous section for $\pi$. Since

$$\pi(s(c^{-1}\alpha(g)^{g}c))=\pi(s(c)^{-1}s\alpha(g)^{g}s(c)),$$ then,  there exists a unique $\gamma(g)\in ker\pi=A$ such that $$\gamma(g)(s(c)^{-1}s\alpha(g)^{g}s(c))=s(c^{-1}\alpha(g)^{g}c).$$ It is clear that the map $\gamma:G\rightarrow A$, $g\mapsto \gamma(g)$ is  continuous. Therefore,\\
\\$\delta^{1}(\beta)(g,h)=s\beta(g).^{g}s\beta(h).(s\beta(gh))^{-1}$ \\ \\$=s(c^{-1}\alpha(g)^{g}c).^{g}s(c^{-1}\alpha(h)^{h}c).(s(c^{-1}\alpha(gh)^{gh}c))^{-1}$
$$=\gamma(g)[s(c)^{-1}s\alpha(g)^{g}s(c)].^{g}(\gamma(h)[s(c)^{-1}s\alpha(h)^{h}s(c)]).(\gamma(gh)[s(c)^{-1}s\alpha(gh)^{gh}s(c)])^{-1}$$ $=\ ^{g}\gamma(h)\gamma(gh)^{-1}\gamma(g)
[s(c)^{-1}s\alpha(g)^{g}s(c)].^{g}[s(c)^{-1}s\alpha(h)^{h}s(c)].[s(c)^{-1}s\alpha(gh)^{gh}s(c)]^{-1}$\\
\\$=\tilde{\delta}^{1}(\gamma)(g,h)[s(c)^{-1}s\alpha(g)^{g}s\alpha(h)(a\alpha(gh))^{-1}s(c)]$\\
\\$=\tilde{\delta}^{1}(\gamma)(g,h)[s(c)^{-1}\delta^{1}(\alpha)(g,h)s(c)]=\tilde{\delta}^{1}(\gamma)(g,h)[\delta^{1}(\alpha)(g,h)].$\\

The last equality is obtained from the fact that
$\delta^{1}(\alpha)(g,h) \in A\subset Z(B)$ and $s(c)\in B$. Now,
note that   $\delta^{1}(\alpha)$ is cohomologous to
$\delta^{1}(\beta)$, when $\alpha$ is cohomologous to $\beta$.
 Thus, $\delta^{1}$ is well-defined.
\section{A Cohomology Exact Sequence}
Let $(X,x_{0}), (Y,y_{0})$ be pointed sets in $\mathcal{PS}$ and  $f:(X,x_{0})\rightarrow (Y,y_{0})$  a pointed map, i.e., $f:X\rightarrow Y$ is a map such that $f(x_{0})=y_{0}$. For simplicity, we write $f:X\rightarrow Y$ instead of $f:(X,x_{0})\rightarrow (Y,y_{0})$. The  kernel of $f$,  denoted by $Ker(f)$, is the set of all points of $X$ that are mapped to the basepoint $y_{0}$.  A sequence $\xymatrix{(X,x_{0}) \ar[r]^{f} & (Y,y_{0})\ar[r]^{g} & (Z,z_{0})}$ of pointed sets and pointed maps is called an exact sequence if $Ker(g)=Im(f)$.

\begin{Theorem}(i) Let $\xymatrix{1\ar[r] & A\ar[r]^{\iota} & B\ar[r]^{\pi} & C\ar[r] & 1}$ be a short exact sequence of topological $G$-modules and continuous $G$-homomorphisms, where $\iota$ is  homeomorphic embedding. Then, the following  is an exact sequence of pointed sets,
$$\xymatrix{0\ar @{->}[r]& H^{0}(G,A)\ar @{->}[r]^{\iota_{0}^{*}}&H^{0}(G,B)\ar @{->}[r]^{\pi_{0}^{*}}&H^{0}(G,C) & & \\ \ar @{->}[r]^{\delta^{0}}&H^{1}(G,A)\ar @{->}[r]^{\iota_{1}^{*}}&H^{1}(G,B)\ar @{->}[r]^{\pi_{1}^{*}}&H^{1}(G,C). }\eqno{(4.1)}$$
(ii) In addition, if $\iota(A)\subset Z(B)$, and $\pi$ has a continuous section, then
$$\xymatrix{ 0\ar @{->}[r]& H^{0}(G,A)\ar @{->}[r]^{\iota_{0}^{*}}&H^{0}(G,B)\ar @{->}[r]^{\pi_{0}^{*}}&H^{0}(G,C)  \\ \ar @{->}[r]^{\delta^{0}}&H^{1}(G,A)\ar @{->}[r]^{\iota_{1}^{*}}&H^{1}(G,B)\ar @{->}[r]^{\pi_{1}^{*}}&H^{1}(G,C)\ar [r]^{\delta^{1}}& H^{2}(G,A)} \eqno{(4.2)}$$
    is an exact sequence of pointed sets.
\end{Theorem}
\textbf{Proof.} (i):  We prove  the exactness term by term.
\par 1. Exactness at $H^{0}(G,A)$: This is clear, since $\iota$ is one to one.
\par 2. Exactness at $H^{0}(G,B)$: Since $\pi_{0}^{*}\iota_{0}^{*}=(\pi\iota)_{0}^{*}=1$, then, $Im(\iota_{0}^{*})\subset Ker(\pi_{0}^{*})$. Now we show that $Ker\pi_{0}^{*}\subset Im\iota_{0}^{*}$. If $b\in Ker\pi^{*}_{0}$, then, $\pi(b)=1$ and $b\in H^{0}(G,B)$. There is an $a\in A$ such that $\iota(a)=b$. Moreover, $\iota(^{g}a)=^{g}\iota(a)=\iota(a)$, $\forall g\in G$. So, $^{g}a=a$, $\forall g\in G$, since $\iota$ is one to one. Thus, $a\in H^{0}(G,A)$. Hence, $b\in Im(\iota^{*}_{0})$.

\par 3. Exactness at $H^{0}(G,C)$: Take $c\in Im(\pi_{0}^{*})$. So, $c=\pi(b)$ for some $b\in H^{0}(G,B)$. Thus, $\delta^{0}(c)(g)=b^{-1}\-^{g}b$. Hence, $\delta^{0}(c)\sim \alpha_{0}$.
Conversely, if $\delta^{0}(c)\sim \alpha_{0}$, then, there is  $a_{1}\in A$ such that $\delta^{0}(c)(g)=a_{1}^{-1}\-^{g}a_{1}$, $\forall g\in G$. Let $c=\pi(b)$ for some $b\in B$. Then, by definition of $\delta^{0}(c)(g)$, there is  $a_{2}\in A$ such that $b^{-1}$$^{g}b=a_{2}^{-1}\delta^{0}(c)(g)^{g}a_{2}$, $\forall g\in G$. So, $b(a_{1}a_{2})^{-1}\in H^{0}(G,B)$.
Since $\pi_{0}^{*}(b(a_{1}a_{2})^{-1})=c$, then, $c\in Im\pi^{*}_{0}$.
\par 4. Exactness at $H^{1}(G,A)$: Let $c\in H^{0}(G,C)$. Then, there is $b\in B$ such that $\pi(b)=c$. So, $$\iota_{1}^{*}\delta^{0}(c)(g)=\iota(\delta^{0}(c)(g))=\iota(b^{-1}\-^{g}b)=b^{-1}\-^{g}b.$$ Consequently, $\iota_{1}^{*}\delta^{0}(c)\sim \beta_{0} $, where $\beta_{0}(g)=1, \forall g\in G$. Conversely, let $\alpha\in Ker\iota_{1}^{*}$. Then, there is  $b\in B$ such that $\iota\alpha(g)=b^{-1}\-^{g}b, \forall g\in G$. So, $\pi(b^{-1}$$^{g}b)=1, \forall g\in G$. Take $c=\pi(b)$. Hence, $c\in H^{0}(G,C)$. Thus, $\delta^{0}(c)\sim \iota(\alpha)=\alpha$.
\par 5. Exactness at $H^{1}(G,B)$: Since $\pi_{1}^{*}\iota_{1}^{*}=(\pi\iota)_{1}^{*}=1$, then, $Im\iota_{1}^{*}\subset Ker\pi_{1}^{*}$. Conversely, let $\beta\in ker\pi_{1}^{*}$. Then, there is  $c\in C$ such that $\pi\beta(g)=c^{-1}$$^{g}c$, for all $g\in G$. Let $b \in B$ and $c=\pi(b)$. Therefore, $\pi(\beta(g))=\pi(b^{-1}$$^{g}b)$, $\forall g\in G$. On the other hand, the map $\tau:A\rightarrow A$, $a\mapsto b^{-1}ab$, is a topological isomorphism, because $A$ is a normal subgroup of $B$. So, for every $g\in G$ there is a
 unique element $a_{g}\in G$ such that, $\beta(g)=(b^{-1}a_{g}b)(b^{-1}$$^{g}b)$.
Thus, $\beta(g)=b^{-1}a_{g}\-^{g}b$, $\forall g\in G$.
 Hence, $a_{g}=b\beta(g)^{g}b^{-1}$, $\forall g\in G$.
Obviously, the map $\alpha:G\rightarrow A$ via $\alpha(g)=a_{g}$ is a continuous derivation,
and $\iota_{1}^{*}(\alpha)\sim \iota_{1}^{*}(\beta)=\beta$.
\par (ii): It is enough to show the exactness at $H^{1}(G,C)$. Let $[\beta] \in H^{1}(G,B)$ and $s$ be a continuous section for $\pi$. Then,  there is a continuous map $z:G\to A$ such that $s\pi\beta(g)=\beta(g)z(g)$. Thus,
\begin{center}
$\delta^{1}(\pi_{1}^{*}(\beta))(g,h)=s(\pi\beta(g))^{g}s(\pi\beta(h))(s(\pi\beta(gh)))^{-1}=\beta(g)^{g}\beta(h)\beta(gh)^{-1}\tilde\delta^{1}(z)(g,h)=\tilde\delta^{1}(z)(g,h)$.
\end{center}
 So, $Im\pi_{1}^{*}\subset Ker\delta^{1}$. Conversely, let $[\gamma]\in ker\delta^{1}$. Then, there is a continuous function $\alpha\in \tilde{C}^{1}(G,A)$ such that $\delta^{1}(\gamma)=\tilde{\delta}^{1}(\alpha)$. Thus,
 \begin{center}
 $s\gamma(g)^{g}s\gamma(h)(s\gamma(gh))^{-1}=\-^{g}\alpha(h)\alpha(gh)^{-1}\alpha(g), \forall g,h\in G$.
 \end{center}
 Assume $\beta(g)=s\gamma(g)\alpha(g)^{-1}, \forall g\in G$. Since $A\subset Z(B)$,
then, $\beta$ is a continuous derivation from $G$ to $B$.  Also $\pi\beta=\gamma$.
 Hence, $\pi_{1}^{*}([\beta])=[\gamma]$.\\
\par The following two corollaries are immediate consequences of Theorem 4.1.
\begin{Corollary} Let $\xymatrix{1\ar[r]  &A\ar[r]^{\iota}& B\ar[r]^{\pi}& C\ar[r]&1}$ be a short exact sequence of discrete $G$-modules, and $G$-homomorphisms then, there is the exact sequence (i) of pointed sets.
\end{Corollary}
\begin{Corollary}Let $\xymatrix{0\ar[r]  &A\ar[r]^{\iota}& B\ar[r]^{\pi}& C\ar[r]&1}$ be a central short exact sequence of discrete $G$-modules and $G$-homomorphisms then, there is the exact sequence (ii) of pointed sets.
\end{Corollary}
\begin{Remark} If we restrict ourselves to the discrete coefficients then, Corollary 4.2 and Corollary 4.3 are the same as Proposition 36 and Proposition 43 in [6, Chapter I], respectively.
\end{Remark}
 \begin{Lemma}Let $G$ be a connected group, and $A$ a totally disconnected abelian topological $G$-module. Then, $H^{n}(G,A)=0$ for every $n\geq1$.
\end{Lemma}
\textbf{Proof.} Consider the coboundary maps $\tilde{\delta}^{n}:\tilde{C}^{n}(G,A)\rightarrow \tilde{C}^{n+1}(G,A)$.
Since $G$ is connected and $A$ is totally disconnected then, $G$ acts trivially on $A$, and the  continuous maps from $G^{n}$ into $A$ are constant. If $n$ is an even positive integer then, one can see that $Ker\tilde{\delta}^{n}=\tilde{C}^{n}(G,A)$ and $Im\tilde{\delta}^{n-1}=\tilde{C}^{n}(G,A)$. Thus, $H^{n}(G,A)=\frac{Ker\tilde{\delta}^{n}}{Im\tilde{\delta}^{n-1}}=\frac{\tilde{C}^{n}(G,A)}{\tilde{C}^{n}(G,A)}=0$. Now suppose that $n$ is odd. It is easy to check that  $Ker\tilde{\delta}^{n}=0$. Consequently, $H^{n}(G,A)=0$.
\begin{Remark} The existence of continuous section in  theorem 4.1 is  essential.
\end{Remark} For example, consider the central short exact sequence of trivial $\Bbb S^{1}$-modules:
\begin{center}
    $\xymatrix{0\ar[r]&\Bbb Z\ar @{->} [r]^{\imath}&\Bbb R\ar [r]^{\pi} &\Bbb S^{1}\ar [r] & 1},$
\end{center}
here $\pi$ is the exponential map, given by $\pi(t)=e^{2\pi it}$ and $\imath$ is the inclusion map. This central exact sequence has no continuous section. For if it has a continuous section then by [1, Lemma 3.5]  , $\Bbb R$ is homeomorphic to $\Bbb Z\times \Bbb S^{1}$. This is a contradiction since $\mathbb{R}$ is connected but $\Bbb Z\times \Bbb S^{1}$ is disconnected . Thus, $Hom_{c}(\mathbb{S}^{1},\mathbb{R})=0$. Now by Lemma 4.5,
\begin{center}
    $H^{1}(\Bbb S^{1},\Bbb Z)=H^{1}(\Bbb S^{1},\Bbb R)=H^{2}(\Bbb S^{1},\Bbb Z)=0$,
\end{center}
On the other hand,
\begin{center}
    $H^{1}(\Bbb S^{1},\Bbb S^{1})=Hom_{c}(\Bbb S^{1},\Bbb S^{1})\neq 0$.
\end{center}
Thus, we don't obtain the exact  sequence (4.2).

\section{Connectivity of Topological Groups}
In this section by using the inflation and the restriction maps, we find a necessary and sufficient condition for connectivity of a topological group $G$.
\begin{Definition} Let $A$ be a topological $G$-module and $A'$ a topological $G'$-module. Suppose that $\phi:G'\rightarrow G$, $\psi:A\rightarrow A'$ are  continuous homomorphisms. Then, we call $(\phi,\psi)$ a cocompatible pair if
$$\psi(^{\phi(g')}a)=\-^{g'}\psi(a), \forall g'\in G', \forall a\in A.$$
\end{Definition}
For example, if $N$ is a subgroup of $G$ and $A$  a topological $G$-module then, $(\imath,Id_{A})$ is a cocompatible pair, where  $\imath:N\rightarrow G$ is the inclusion map and $Id_{A}$ is the identity map. Also, suppose that $\pi:G\rightarrow G/N$ is the natural projection  and $\jmath:A^{N}\rightarrow A$ is the inclusion map. Then, $(\pi,\jmath)$ is a cocompatible pair.
\\
\par Note that a cocompatible pair $(\phi,\psi)$ induces  a natural map as follows:
\begin{center}
    $Der_{c}(G,A)\rightarrow Der_{c}(G',A')$ by  $\alpha\mapsto \psi\alpha\phi$,
\end{center}
which induces the map:

\begin{center}
     $(\phi,\psi)^{*}:H^{1}(G,A)\rightarrow H^{1}(G',A')$ by $[\alpha]\mapsto [\psi\alpha\phi]$.
\end{center}
\begin{Definition} Let $N$ be a subgroup of $G$ and $A$  a topological $G$-module. Suppose that $\imath:N\rightarrow G$ is the inclusion map. The induced map $(\imath,Id_{A})^{*}$ is called the \emph{restriction map} and it is denoted by $Res^{1}:H^{1}(G,A)\rightarrow H^{1}(N,A)$.
\end{Definition}
\begin{Definition} Let $N$ be a normal subgroup of $G$ and
 $A$  a topological $G$-module. Suppose that $\pi:G\rightarrow G/N$ is the natural projection  and $\jmath:A^{N}\rightarrow A$ is the inclusion map.
The induced map $(\pi,\jmath)^{*}$ is called the \emph{inflation map} and it is denoted by $Inf^{1}:H^{1}(G/N,A^{N})\rightarrow H^{1}(G,A)$.
\end{Definition}
\par Note that if $A$ is an abelian topological $G$-modules then, $Inf^{1}$ and $Res^{1}$ are group homomorphisms.
%
%
%
%
%
\begin{Lemma} Let $A$ be a topological $G$-module, and $N$ a normal subgroup of $G$. Then,
\begin{itemize}
\item[(i)]  \emph{$H^{1}(N,A)$ is a $G/N$-set. Moreover, if $A$ is an abelian topological $G$-module then, $H^{1}(N,A)$ is an abelian $G/N$-module.}

\item[(ii)] $ImRes^{1}\subset H^{1}(N,A)^{G/N}$.

\end{itemize}
\end{Lemma}
\textbf{Proof.} (i) Since $N$ is a normal subgroup of $G$, then, there is an action of $G$ on $Der_{c}(N,A)$ as follows:\\
 For every $g\in G$ we define $^{g}\alpha=\tilde{\alpha}, \forall g\in G$, with $\tilde{\alpha}(n)=\-^{g}\alpha(^{g^{-1}}n), n\in N$.
\par In fact, $\tilde{\alpha}$ is continuous and we have:
\begin{center}
    $\tilde{\alpha}(mn)=\-^{g}\alpha(^{g^{-1}}(mn))=\-^{g}\alpha(^{g^{-1}}m^{g^{-1}}n)=
    \-^{g}\alpha(^{g^{-1}}m)\-^{mg}\alpha(^{g^{-1}}n)=\tilde{\alpha}(m)^{m}\tilde{\alpha}(n)$,
\end{center}
whence, $\tilde{\alpha}\in Der_{c}(N,A)$. It is clear that $^{gh}\alpha=\-^{g}(^{h}\alpha)$. Moreover, if $A$ is an abelian group, it is easy to verify that  $^{g}(\alpha\beta)=\-^{g}\alpha^{g}\beta$. Now suppose that $\alpha\sim\beta$. Then, there is an $a\in A$ with $\beta(n)=a^{-1}\alpha(n)^{n}a, \forall n\in N$. Thus, for every $g\in G$, $n\in N$,
\begin{center}
 $^{g}\beta(^{g^{-1}}n)=\-^{g}a^{-1}(^{g}\alpha(^{g^{-1}}n))^{g}(^{^{g^{-1}}n}a)$.
\end{center}
Therefore,
\begin{center}
    $\tilde{\beta}(n)=(^{g}a)^{-1}\tilde{\alpha}(n)^{n}(^{g}a)$, i.e., $\tilde{\alpha}\sim \tilde{\beta}$.
\end{center}
Thus, the action of $G$ on $Der_{c}(N,A)$ induces an action of $G$ on $H^{1}(N,A)$. It is sufficient to show for every $m\in N$, $^{m}\alpha\sim \alpha$. In fact, for every $n\in N$
\begin{center}
    $^{m}\alpha(^{m^{-1}}n)=\-^{m}\alpha(m^{-1}nm)=\ ^{m}(\alpha(m^{-1})^{m^{-1}}\alpha(n)^{m^{-1}n}\alpha(m))
    =\ ^{m}\alpha(m^{-1})\alpha(n)^{n}\alpha(m)=\alpha(m)^{-1}\alpha(n)^{n}\alpha(m)$.
\end{center}
Thus, $\tilde{\alpha}\sim \alpha$.
\par (ii) By a similar argument as in (i), we have
\begin{center}
    $^{g}\alpha(^{g^{-1}}n)=\alpha(g)^{-1}\alpha(n)^{n}\alpha(g)$, $\forall g\in G, n\in N$
\end{center}
whence, $^{gN}(\alpha\imath)\sim \alpha\imath, \forall gN\in G/N$.
%
%
%
%
%
%
%
%
%
%
%
%
%
%
\begin{Lemma} Let $N$ be a normal subgroup of a topological group $G$ and $A$  a topological $G$-module. Then,  there is an exact sequence
\begin{center}
    $\xymatrix{1\ar[r]&H^{1}(G/N,A^{N})\ar[r]^{Inf^{1}}&H^{1}(G,A)\ar[r]^{Res^{1}}& H^{1}(N,A)^{G/N}}$.
\end{center}
\end{Lemma}
\textbf{Proof.} The map $Inf^{1}$ is one to one: If $\alpha, \beta\in Der_{c}(G/N,A^{N})$ and $Inf^{1}[\alpha]=Inf^{1}[\beta]$, then, $\alpha\pi\sim \beta\pi$. Thus, there is an $a\in A$ such that $\beta\pi(g)=a^{-1}\alpha\pi(g)^{g}a, \forall g\in G$. Hence, $\beta(gN)=a^{-1}\alpha(gN)^{g}a,  \forall gN\in G/N$.On the other hand, if $g\in G$, then, $\alpha(gN)=\beta(gN)=1$, and hence, $a\in A^{N}$. This implies that $^{(gN)}a=\-^{g}a, \forall g\in G$. Consequently, $\alpha \sim \beta$, i.e., $Inf^{1}$ is one to one.
\par Now we show that  $KerRes^{1}= ImInf^{1}$. Since $Res^{1}Inf^{1}[\alpha]=[\alpha(\pi\imath)]=1$, then, $ImInf^{1}\subset KerRes^{1}$.
\par Let $[\alpha] \in KerRes^{1}$. Then, there is an $a\in A$ such that $\alpha(n)=a^{-1}$$^{n}a$, $\forall n\in N$. Consider the continuous derivation $\beta$ with $\beta(g)=a\alpha^{g}a^{-1}$, $\forall g\in G$. Since $\beta(n)=1, \forall n\in N$ then, $\beta$  induces  the continuous derivation $\gamma:G/N\rightarrow A$ via $\gamma(gN)=\beta(g)$. Also $Im\gamma\subset A^{N}$, since for all $n\in N$,
\begin{center}
    $^{n}\gamma(gN)=\-^{n}\beta(g)=\beta(ng)=\beta(g)^{g}\beta(g^{-1}ng)=\beta(g)=\gamma(gN)$.
\end{center}
Hence, $Inf^{1}[\gamma]=[\gamma\pi]=[\beta]=[\alpha]$. Consequently, $KerRes^{1}\subset ImInf^{1}$.
\begin{Lemma} Let $G$ be a topological group and $A$ a topological $G$-module. Suppose that $A$ is totally disconnected and $G_{0}$ the identity component of $G$. Then,  the map
$$\xymatrix{H^{1}(G/G_{0},A)\ar[r]^{Inf^{1}} & H^{1}(G,A)}$$
is bijective.
\end{Lemma}
 \textbf{Proof.} Since $G_{0}$  acts trivially on $A$, then,  $A^{G_{0}}=A$. On the other hand, $H^{1}(G_{0},A)=1$. Thus, by Lemma 5.5, the sequence
$$\xymatrix{0\ar[r] & H^{1}(G/G_{0},A)\ar[r]^{Inf^{1}} & H^{1}(G,A)\ar[r] & 0}$$
is exact.
\begin{Theorem} Let $G$ be a topological group which  has an open component. Then, $G$ is connected iff $H^{1}(G,A)=1$ for every discrete abelian $G$-module $A$.
\end{Theorem}
\textbf{Proof.} Assume $G$ is a connected group and $A$ a discrete abelian $G$-module. Since every discrete $G$-module $A$ is totally disconnected then, $H^{1}(G,A)=1$.
 \par Conversely, Suppose that $H^{1}(G,A)=1$, for every discrete abelian $G$-module $A$.  By Lemma 5.6,  $H^{1}(G/G_{0},A)=1$,  for every discrete abelian  $G$-module $A$. Since $G/G_{0}$ is discrete, then, the cohomological dimension of  $G/G_{0}$ is equal to  0 which implies that  $G/G_{0}=1$ [4, Chapter VIII],
  i.e., $G=G_{0}$.
\section{Complements and First Coholomology}
\par Let $G$ and $A$ be topological groups. Suppose that $\chi:G\times A\rightarrow A$ is a continuous
map such that  $\tau_{g}:A\rightarrow A$, defined by
$\tau_{g}(a)=\chi(g,a)$, is a  homeomorphic automorphism of $A$ and the map $g\mapsto \tau_{g}$ is a
homomorphism of $G$ into the group of homeomorphic automorphisms, $Aut_{h}(A)$, of $A$.
By $G\ltimes_{\chi} A$  we mean the (topological) semidirect product with the group
operation, $(g, a)(h, b) = (gh, \tau_{h}(a)b)$, and  the product topology of $G\times A$. Sometimes for simplicity we denote $G\ltimes_{\chi} A$ by $G\ltimes A$  and  view $G$ and $A$ as
topological subgroups of $G\ltimes A$ in a natural way. So every element $e$ in $G\ltimes N$
can be written uniquely as $e=gn$ for some $g\in G$ and $n\in N$.
\\
\par Let $E=G\ltimes N$.
A subgroup $X$ of $E$ such that $E\simeq X\ltimes N$ is called a complement of $N$ in $E$. Indeed, any conjugate of $G$ is a complement.
\par We  show that the complements of $N$ in $E$ correspond to continuous derivations from $G$ to $N$. If $X$ is any complement, for every $g\in G$, then, $g^{-1}$ has a unique expression  of the form $g^{-1}=xn$ where $x\in X$ and $n\in N$. Define $\alpha_{X}:G\rightarrow N$ by $\alpha_{X}(g)=n$. Obviously, $\alpha_{X}(g)=\pi_{2}|_{G}(g^{-1})$, where $\pi_{2}:X\ltimes N\rightarrow N$ is given by $\pi_{2}(x,n)=n$. Hence, $\alpha_{X}$ is  continuous. Now if $g_{i}\in G$ then, $g_{i}^{-1}=x_{i}n_{i}$ for some $x_{i}\in X$ , $n_{i}\in N$, $i=1, 2$. We have:
$$(g_{1}g_{2})^{-1}=g_{2}^{-1}g_{1}^{-1}=x_{2}n_{2}x_{1}n_{1}=x_{2}x_{1}^{x_{1}^{-1}}n_{2}n_{1}=x_{2}x_{1}n_{1}^{(n_{1}^{-1}x_{1}^{-1})}n_{2}=x_{2}x_{1}n_{1}^{g_{1}}n_{2}.$$
By definition of $\alpha_{X}$,  $\alpha_{X}(g_{1}g_{2})=\alpha_{X}(g_{1})^{g_{1}}\alpha(g_{2})$, i.e., $\alpha_{X}\in Der_{c}(G,N)$.
\par So, we have associated a continuous derivation with each complement. Conversely, suppose that $\alpha:G\rightarrow N$ is a  continuous derivation. Then, $X_{\alpha}=\{\alpha(g)g| g\in G\}\subset E$  is a  corresponding complement to $\alpha$ in $E$. Obviously, the continuous map $\kappa: g\mapsto \alpha(g)g$, is a homomorphism. Suppose that $\pi_{1}:G\ltimes N\rightarrow G$ is given by $\pi_{1}(g,n)=g$. Hence, $\pi_{1}|_{X_{\alpha}}:X_{\alpha}\rightarrow G$ is the inverse of $\kappa$, since $\pi_{1}|_{X_{\alpha}}(\alpha(g)g)=\pi_{1}|_{X_{\alpha}}(g^{g^{-1}}n)=g, \forall g\in G$. Thus, $X_{\alpha}\simeq G$.
\par Define the map $\chi:X_{\alpha}\times N \rightarrow N$ by $\chi(\alpha(g)g)=\-^{g}n$, for all $g\in G, n\in N$. Clearly, $\chi$ is a continuous map. Hence, $X_{\alpha}\ltimes_{\chi} N \simeq G\ltimes N=E$.
\\
\par In fact we have proved the following theorem.
 \begin{Theorem}Let $G$ be a topological group and $N$ a topological $G$-module. Then, the map $X\mapsto \alpha_{X}$ is a bijection from the set of all complements of $N$ in $G\ltimes N$ onto $Der_{c}(G,N)$.
\end{Theorem}
\begin{Theorem}If $A$ is a topological $G$-module then, there is a map from $H^{1}(G,A)$ onto the set of conjugacy classes of complements of $A$ in $G\ltimes A$. Moreover, if $A$ is an abelian group then, this map is one to one.
\end{Theorem}
\textbf{Proof.} Suppose that $X$ and $Y$ are the complements of $A$ in $G\ltimes A$ such that  $\alpha_{X}\sim \alpha_{Y}$. Hence, there is  $a\in A$ such that $\alpha_{Y}(g)=a^{-1}\alpha_{X}(g)^{g}a$, $\forall g \in G$. Thus, for each $g\in G$, we have $\alpha_{Y}(g)g=a^{-1}\alpha_{X}(g)^{g}ag=a^{-1}\alpha_{X}(g)ga$. This implies that $X=\-^{a^{-1}}Y$.
\par Moreover, suppose that $A$ is an abelian group and  $X$ and $Y$ are conjugate complements. So, $X=\-^{n}Y$ for some $n\in N$. If $g\in G$, then, $\alpha(g)g\in X$ where $\alpha_{X}$ is a continuous derivation arising from $X$. Hence, $\alpha_{X}(g)g=\-^{n}y$ for some $y\in Y$. Now $^{n}y=[n,y]y$, so, $\alpha_{X}(g)g=[n,y]y$, which shows that

$$[n,g]=ngn^{-1}g^{-1}=n(\alpha_{X}(g)^{n}y)n^{-1}(\alpha_{X}(g)^{n}y)^{-1}=$$ $$n(\alpha_{X}(g))nyn^{-1}n^{-1}ny^{-1}n^{-1}(\alpha_{X}(g))^{-1}
=$$ $$(n\alpha_{X}(g))(nyn^{-1}y^{-1})(n^{-1}\alpha_{X}(g)^{-1})=[n,y]=\-^{n}yy^{-1}$$

because $A$ is an abelian group. Therefore,
\begin{center}
$g^{-1}=y^{-1}[n,y]^{-1}\alpha_{X}(g)=y^{-1}(y^{n}y^{-1}\alpha_{X}(g))=y^{-1}(y\alpha_{X}(g)^{n}y^{-1})$.
\end{center}
Thus, by definition of $\alpha_{Y}$, we get $\alpha_{Y}(g)=y\alpha_{X}(g)^{n}y^{-1}$.
 Consequently, $\alpha_{X}\sim\alpha_{Y}$.
\\
\par As an immediate result, we have the following corollary.
\begin{Corollary} Let $A$ be a topological $G$-module and $H^{1}(G,A)=1$. Then,
the complements of $A$ in $G\ltimes A$ are conjugate.
\end{Corollary}
\section{Vanishing of $H^{1}(G,A)$}
Let $G$ be a compact Hausdorff group and  $A$  a topological $G$-module. Suppose that $A$ is an  almost connected locally compact Hausdorff group. Then, we prove  there exists a $G$-invariant maximal compact subgroup $K$ of $A$, and for every such  topological submodule $K$, the natural map $\iota_{1}^{*}:H^{1}(G,K)\rightarrow H^{1}(G,A)$ is onto. In addition, as a result, If  $A$ has trivial maximal compact subgroup then, $H^{1}(G,A)=1$.\\
Recall that $G$ is almost connected if $G/G_{0}$ is compact where $G_{0}$ is the connected component of the identity of $G$.
\begin{Definition}An element $g\in G $ is called periodic if it is contained in a compact subgroup of $G$. The set of all periodic elements of $G$ is denoted by $P(G)$.
\end{Definition}
 \begin{Definition} A maximal compact subgroup $K$ of a topological group $G$ is a subgroup $K$ that is a compact space in the subspace topology, and maximal amongst such subgroups.
\end{Definition}
\par If a topological group $G$ has a maximal compact subgroup $K$, then, clearly $gKg^{-1}$ is a maximal compact subgroup of $G$ for any $g\in G$. There exist topological groups with maximal compact subgroups and
compact subgroups which are not contained in any maximal one [3]. Note that if $G$ is almost connected then, $P(G/G_{0})=G/G_{0}$.
\begin{Lemma} Let $G$ be a locally compact topological group such that $P(G/G_{0})$ is a compact subgroup of $G/G_{0}$, and  $K$  a maximal compact subgroup of $G$. Then, any compact subgroup of $G$ can be conjugated into $K$ [3, Theorem 1].
\end{Lemma}
\begin{Lemma}Let $G$ be a compact group and $A$  a topological $G$-module such that $A$ is a locally compact almost connected, and let $C$ be a $G$-invariant compact subgroup of $A$. Then, there exists a $G$-invariant maximal compact subgroup $K$ of $A$ which contains $C$.
\end{Lemma}
 \textbf{Proof.} Let $E = G\ltimes A$, be the semidirect product of $A$ and $G$ with respect to the action of $G$ on $A$. Note that topologically $E$ is the product of $A$ and $G$. We first observe that $E/E_{0}$ is almost connected. Let $A_{0}$, $G_{0}$ and $E_{0}$ be the components of $A, G$ and $E$, respectively. It is easily seen that $E_{0} = A_{0}\times G_{0}$. Also $E/(A_{0}\times G_{0})$ is  homeomorphic to the compact space $A/A_{0}\times G/G_{0}$. Hence,  $E/E_{0}$ is compact. Consequently, $E$ is almost connected. Now, by assumption, $C$ is a $G$-invariant compact subgroup of $A$. Thus, $G\ltimes C$ is a compact subgroup of $E$. Since $E$ is almost connected, there exists a maximal compact subgroup $L$ of $E$ which contains $G\ltimes C$. Let $K=L\cap A$. Since $K$ is a closed subspace of $L$, then, $K$ is compact.  Also $L$ contains $G$. Thus, $L$ is $G$-invariant. In fact, for every $g \in G$ and every $\ell \in L$, we have $^{g}\ell=g\ell g^{-1} \in L$. This immediately implies that $K$ is $G$-invariant, since $A$ is $G$-invariant.  Let $K'$ be a compact subgroup of $G$. By Lemma 7.3, there is $e\in E$ such that $eK'e^{-1}\subset L$. Thus, $eK'e^{-1}\subset L\cap A=K$. But there exist $g\in G$ and $ a\in A$ such that $e=ga$. Thus, $aK'a^{-1}\subset g^{-1}Kg=\-^{g^{-1}}K=K$. Therefore,  $K$ is a $G$-invariant maximal compact subgroup of $A$ which contains $C$.
\begin{Theorem} Let $G$ be a compact Hausdorff group and  $A$  a topological $G$-module. Let  $A$ be an  almost connected locally compact Hausdorff group. Then, there exists a $G$-invariant maximal compact subgroup $K$ of $A$, and for every such  topological submodule $K$, the natural map $\iota_{1}^{*}:H^{1}(G,K)\rightarrow H^{1}(G,A)$ is onto.
\end{Theorem}
 \textbf{Proof.} By Lemma 7.4, there exists a $G$-invariant maximal compact subgroup $K$ of $A$. Also $G\ltimes K$ is a maximal compact subgroup of $G\ltimes A$ [2, Theorem 1.1]. Let $\alpha:G\rightarrow A$ be a continuous derivation. Then, define the continuous homomorphism $\kappa:G\rightarrow G\ltimes A$ via $g\mapsto \alpha(g)g$. Since $\kappa$ is a continuous homomorphism then, $\kappa(G)$ is a compact subgroup of $G\ltimes A$. By Lemma 7.3 there is $ag\in G\ltimes A$ such that $(ag)\kappa(G)(ag)^{-1}\subset G\ltimes K, \forall x \in G$. This is equivalent to $(ag)\alpha(x)x(ag)^{-1}\subset G\ltimes K, \forall x \in G$. Hence, for all $x\in G$, $^{g}[^{g^{-1}}a\alpha(x)^{x}(^{g^{-1}}a^{-1})]gxg^{-1} \in G\ltimes K$. Since $K$ is $G$-invariant then, $(^{g^{-1}}a)\alpha(x)^{x}(^{g^{-1}}a^{-1})\in K, \forall x\in G$. Now define $\beta:G\rightarrow K$ by $\beta(x)=(^{g^{-1}}a)\alpha(x)^{x}(^{g^{-1}}a^{-1}), \forall x\in G $. Hence, $\iota_{1}^{*}([\beta])=[\alpha]$, i.e., $\iota_{1}^{*}$ is onto map.
 \begin{Corollary} Let $G$ be a compact Hausdorff group and  $A$ a topological $G$-module. Let  $A$ be  an  almost connected locally compact Hausdorff group with the trivial maximal compact subgroup. Then, $H^{1}(G,A)=1$.
\end{Corollary}
 \textbf{Proof.} It is clear.

\end{document}